\documentclass[11pt]{article}
\usepackage{geometry}
\usepackage{mathrsfs}
\usepackage{color}
\usepackage{epsfig}
\usepackage{graphicx}
\usepackage{epstopdf}
\usepackage{cite}
\usepackage{mathrsfs}
\usepackage{amsfonts}
\usepackage{amsmath}
\usepackage{amsfonts,amssymb,color}
\usepackage{dsfont}
\usepackage{curves}
\usepackage{mathrsfs}
\usepackage{pifont}
\usepackage{amssymb}
\usepackage{latexsym,amsmath,amssymb,amsfonts,epsfig,graphicx,cite,psfrag}
\usepackage{eepic,color,colordvi,amscd}
\def\qed{\hfill \rule{4pt}{7pt}}
\def\pf{\noindent {\it Proof. }}

\renewcommand{\baselinestretch}{1.2}

\title{On the chromatic number of some $P_5$-free graphs\footnote{Supported by NSFC No. 11931106, }}
\author{Wei Dong$^{1,}$\footnote{Email: weidong@njxzc.edu.cn }, Baogang Xu$^{2,}$\footnote{Email: baogxu@njnu.edu.cn, or baogxu@hotmail.com.} and Yian Xu$^{3,}$\footnote{Email: yian$\_$xu@seu.edu.cn}
\\ $^{1}$\small School of Information and Engineering
\\ \small Nanjing Xiaozhuang University, Nanjing,
211171, China
\\ $^{2}$\small Institute of Mathematics, School of Mathematical Sciences
\\ \small Nanjing Normal University, Nanjing, 210023, China
\\$^{3}$\small School of Mathematics, Southeast University, 2 SEU Road, Nanjing, 211189, China}

\date{}
\newtheorem{theorem}{Theorem}[section]

\newtheorem{lemma}{Lemma}[section]

\newtheorem{prob}{Problem}[section]
\newtheorem{conj}{Conjecture}[section]

\newtheorem{obs}{Observation}[section]
\def\qed{\hfill \rule{4pt}{7pt}}
\renewcommand{\baselinestretch}{1.15}
\begin{document}
\maketitle
\begin{abstract}
Let $G$ be a graph. We say that $G$ is perfectly divisible if for each induced subgraph $H$ of $G$, $V(H)$ can be partitioned into $A$ and $B$ such that $H[A]$ is perfect and $\omega(H[B])<\omega(H)$. We use $P_t$ and $C_t$ to denote a path and a cycle on $t$ vertices, respectively.  For two disjoint graphs $F_1$ and $F_2$, we use $F_1\cup F_2$ to denote the graph with vertex set $V(F_1)\cup V(F_2)$ and edge set $E(F_1)\cup E(F_2)$, and use $F_1+F_2$ to denote the graph with vertex set $V(F_1)\cup V(F_2)$ and edge set $E(F_1)\cup E(F_2)\cup \{xy\;|\; x\in V(F_1)\mbox{ and } y\in V(F_2)\}$.  In this paper, we prove that (i) $(P_5, C_5, K_{2, 3})$-free graphs are perfectly divisible, (ii) $\chi(G)\le 2\omega^2(G)-\omega(G)-3$ if  $G$ is   $(P_5,  K_{2,3})$-free with $\omega(G)\ge 2$, (iii) $\chi(G)\le {3\over 2}(\omega^2(G)-\omega(G))$ if  $G$ is   $(P_5, K_1+2K_2)$-free, and (iv) $\chi(G)\le 3\omega(G)+11$ if  $G$ is   $(P_5, K_1+(K_1\cup K_3))$-free.
\begin{flushleft}
\end{flushleft}
\end{abstract}

\section{Introduction}

All graphs considered in this paper are finite, simple, and connected.
Let $G$ be a graph. The {\it clique number $\omega(G)$} of $G$ is the maximum size of the cliques of $G$, and the  {\it independent number $\alpha(G)$} of $G$ is the maximum size of the independent sets of $G$. We use $P_k$ and $C_k$ to denote a \textit{path} and a \textit{cycle} on $k$ vertices respectively. The complete bipartite graph with partite sets of size $p$ and $q$ is denoted by $K_{p, q}$, and the complete graph with $l$ vertices is denoted by $K_l$.

Let $G$ and $H$ be two vertex disjoint graphs. The \textit{union} $G\cup H$ is the graph with $V(G\cup H)= V(G)\cup V(H)$
and $E(G \cup H) = E(G)\cup E(H)$. Similarly, the \textit{join} $G + H$ is the graph with $V(G + H) = V(G)\cup V(H)$ and $E(G + H) =
E(G) \cup E(H)\cup \{xy| \mbox{for each pair} \ x \in V(G) \mbox{ and } y \in V(H)\}$. For positive integer $k$, $kG$ denotes the union of $k$ copies of $G$.

We say that $G$ induces $H$ if $G$ has an induced subgraph isomorphic to $H$, and say that $G$ is $H$-{\em free} if $G$ does not induce $H$. Let $\mathcal{H}$ be a family of graphs. We say that $G$ is ${\cal H}$-free if $G$ induces no member of ${\cal H}$. For a subset $X\subseteq V(G)$, let $G[X]$ denote the subgraph of $G$ induced by $X$. A \textit{hole} of $G$ is an induced cycle of length at least 4, and a $k$-{\em hole} is a hole of length $k$. A $k$-{\em hole} is said to be an {\em odd (even) hole} if $k$ is odd (even). An \textit{antihole} is the complement of some hole. An {\em odd} (resp. {\em even}) antihole is defined analogously.

A coloring of $G$ is an assignment of colors to the vertices of $G$ such that no two adjacent
vertices receive the same color. The minimum number of colors required to color $G$ is said to be the {\em chromatic number} of $G$, denoted by $\chi(G)$. Obviously we have that $\chi(G)\ge \omega(G)$. However, determining the upper bound of the chromatic number of some family of graphs $G$, especially, giving a function of $\omega(G)$ to bound $\chi(G)$ is generally very difficult. Throughout the literature, plenty of work has been taken to investigate this problem. A family $\mathcal{G}$ of graphs is said to be $\chi$-{\em bounded} if there is a function $f$ such that
$\chi(G) \le f(\omega(G))$ for every $G\in \mathcal{G}$, and if such a function $f$ does exist to $\mathcal{G}$, then $f$ is said to be a {\em binding function} of $\mathcal{G}$ \cite{gyarfas1}. A graph $G$ is said to be \textit{perfect} if  $\chi(H)=\omega(H)$ for each induced subgraph $H$. Thus the binding function for perfect graphs is $f(x)=x$. The famous {\em Strong Perfect Graph Theorem} \cite{CRST06} states that a graph is perfect if and only if it induces neither an odd hole nor an odd antihole.  Erd\H{o}s \cite{E59} showed that for any positive integers $k$ and $\l$, there exists a graph $G$ with $\chi(G)\ge k$ and no cycles of length less than $\l$. This result motivates the study of the chromatic number of $\mathcal{H}$-free graphs for some $\mathcal{H}$. Gy\'{a}rf\'{a}s \cite{gyarfas1, G87}, and Sumner \cite{sumner1} independently, proposed the following conjecture.

\begin{conj}{\em \cite{G87, sumner1}}\label{Gyarfas87}
For every tree $T$, $T$-free graphs are $\chi$-bounded.
\end{conj}

Gy\'{a}rf\'{a}s \cite{G87} proved that $\chi(G) \leq {(k-1)}^{\omega(G) - 1}$ for $k\ge 4$ if $G$ is $P_k$-free and $\omega(G)\ge 2$.  Gy\'{a}rf\'{a}s also   suggested that there might exist $\chi$-binding function for these classes of graphs with a better magnitude.

Since $P_4$-free graphs are perfect, determining an optimal binding function of $P_5$-free graphs attracts much attention. Sumner \cite{sumner1} showed that all  $(P_5, K_3)$-free graphs are 3-colorable, and there exist many $(P_5, K_3)$-free graphs with chromatic number 3. Up to now, the best known upper bound for $P_5$-free graphs is due to Esperet {\em et al} \cite{ELMM13}, who showed that if $G$ is $P_5$-free and $\omega(G)\ge 3$ then $\chi(G)\le 5\cdot3^{\omega(G)-3}$, and the bound is sharp for $\omega(G)=3$. A natural question is whether the exponential bound can be improved.

\begin{prob}{\em \cite{S16}}\label{Schiermeyer16}
Are there polynomial functions $f_{P_k}$ for $k\ge 5$ such that $\chi(G)\le f_{P_k}(\omega(G))$ for  every $P_k$-free graph $G$$?$
\end{prob}

\begin{conj}{\em \cite{CKS07}}\label{Choudum06}
There exists a constant $c$ such that for every $P_5$-free graph $G$, $\chi(G)\le c\omega^2(G)$.
\end{conj}

We say that a graph $G$ admits a {\em perfect division} $(A, B)$ if $V(G)$ can be partitioned into $A$ and $B$ such that $G[A]$ is perfect and $\omega(G[B])<\omega(G)$. A graph $G$ is said to be \textit{perfectly divisible} if each of its induced subgraphs admits a   perfect division \cite{H18}.   Obviously, if $G$ is perfectly divisible, then $\chi(G)\le \omega(G)+(\omega(G)-1)+\cdots +2+1={\omega(G)+1\choose 2}$.

Plenty of articles around the above topics have been published in the decades. Here we list some results related to $(P_5, H)$-free graphs for some small graph $H$, and refer the readers to \cite{RS04,SS20,SR19} for more information on Conjecture~\ref{Gyarfas87} and related problems.

\begin{figure}[htbp]
\begin{center}
\includegraphics[scale=0.4]{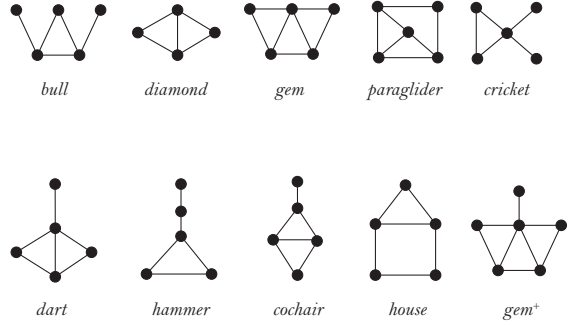}
\end{center}
 {\small \caption{Illustration of some forbidden configurations}}
\label{fig-1}
\end{figure}

A \textit{bull} is a graph consisting of a triangle with two disjoint pendant edges, a {\em cricket} is a graph consisting of a triangle with two adjacent pendant edges, a {\em diamond} is the graph $K_1+P_3$,  a {\em cochair} is the graph obtained from a diamond by adding a pendent edge to a vertex of degree 2,  a {\em dart} is the graph $K_1+(K_1\cup P_3)$, a {\em hammer} is the graph obtained by identifying one vertex of a $K_3$ and one end vertex of a $P_3$, a {\em house} is just the complement of $P_5$,  a {\em gem} is the graph $K_1+P_4$, a {\em gem}$^+$ is the graph $K_1+(K_1\cup P_4)$, and a {\em paraglider} is the graph obtained from a diamond by adding a vertex joining to its two vertices of degree 2 (see Figure~\ref{fig-1}).

Fouquet {\it et al} \cite{FGMT95} proved that ($P_5$, house)-free graphs are perfectly divisible. Schiermeyer \cite{S16} proved that $\chi(G)\le \omega^2(G)$ for ($P_5$, $H$)-free graphs $G$, where $H$ is a graph in $\{$cricket, dart, diamond, gem,  gem$^+$, $K_{1, 3}\}$.
Brause {\em et al} \cite{BRSV13} proved that $\chi(G)\le {\omega(G)+1\choose 2}$ if $G$ is $(P_5$, hammer)-free, Chudnovsky and Sivaraman \cite{CS19} showed that ($P_5$, bull)-free graphs and (odd hole, bull)-free graphs are both perfectly divisible, and Ho\'{a}ng \cite{H18} showed that every (odd holes, banner)-free graph is perfectly divisible. Dong and Xu \cite{DX21+} proved that $(P_5, F)$-free graphs are perfectly divisible, where $F$ is either a cochair or a cricket. Chudnovsky {\it et al} \cite{CKMM20} proved that $\chi(G)\le \lceil{5\omega(G)\over 4}\rceil$ if $G$ is ($P_5$, gem)-free, which improves the results of \cite{CHM19} and \cite{CKS07}.  Char and Karthick \cite{ACTK20} showed that if $G$ is $(P_5$, $K_1+C_4$)-free, then $\chi(G)\le {3\omega(G)\over 2}$.  Huang and Karthick \cite{SHTK19} showed that if $G$ is  $(P_5$, paraglider)-free, then $\chi(G)\le \lceil{3\omega(G)\over 2}\rceil$.

Chudnovsky and Sivaraman \cite{CS19} showed that $\chi(G)\le 2^{\omega(G)-1}$ if $G$ is ($P_5, C_5)$-free, Brause {\em et al} \cite{BDS16} proved that $\chi(G)\le d\cdot\omega^3(G)$ for some constant $d$ if $G$ is $(P_5, K_{2,3})$-free, and Schiermeyer \cite{S17} proved that $\chi(G)\le c\cdot\omega^3(G)$ for some constant $c$ if $G$ is $(P_5, K_1+2K_2)$-free.  In this paper, we study a subclasses of $P_5$-free graphs, and prove the following theorems, which improve some results of \cite{BDS16, S17, W80}.

\renewcommand{\baselinestretch}{1}
\begin{theorem}\label{main0}
Every $(P_5, C_5, K_{2,3})$-free graph is perfectly divisible.
\end{theorem}

\begin{theorem}\label{main1}
If $G$ is $(P_5, K_{2,3})$-free then $\chi(G)\le 2\omega^2(G)-\omega(G)-3$.
\end{theorem}

\begin{theorem}\label{main2}
If $G$ is $(P_5, K_1+2K_2)$-free with $\omega(G)\ge 2$ then $\chi(G)\le {3\over 2}(\omega^2(G)-\omega(G))$.
\end{theorem}

\begin{theorem}\label{main3}
If $G$ is $(P_5, K_1+(K_1\cup K_3))$-free then $\chi(G)\le 3\omega(G)+11$.
\end{theorem}  \renewcommand{\baselinestretch}{1.2}

Theorem~\ref{main1} improves a result of Brause {\em et al} \cite{BDS16} and the upper bound $2\omega^2(G)-\omega(G)-3$ is sharp in the sense that all $(P_5, K_3)$-free graphs are 3-colorable and there are $(P_5, K_3)$-free graphs with chromatic number 3, Theorem~\ref{main2} improves a result of Schiermeyer \cite{S17}, and Theorem~\ref{main3} improves a result of \cite{W80} which states that $\chi(G)\le \frac{1}{2}(\omega^2(G)+\omega(G))$ for $\{2K_2, K_1+(K_1\cup K_3)\}$-free graphs.

It is known (see {\bf Theorem} 14 of \cite{BRSV13}) that the class of $2K_2 \cup 3K_1$-free graphs does not admit a linear binding function, and so one can not expect a linear binding function  for  $(P_5, K_{2,3})$-free graphs or for $(P_5, K_1+2K_2)$-free graphs.

In Section 2, we introduce a few more notations, and list several useful lemmas. Section 3 is devoted to the proof of Theorem~\ref{main0}. Theorems~\ref{main1}, \ref{main2}, and \ref{main3} are proved in Sections 4, 5, and 6 respectively.

\section{Preliminary and Notations}

Let $G$ be a graph, and let $A$ be an  antihole of $G$ with $V(A)=\{v_1, v_2, \cdots, v_h\}$. We always enumerate the vertices of $A$ cyclically  such that $v_iv_{i+1}\notin E(G)$, and simply write $A=v_1v_2\cdots v_h$. In this paper, the summations of subindex are taken modulo $h$ for some $h$, and we always set $h+1\equiv 1$.

\begin{obs}\label{obs1}
The vertices of an odd antihole cannot be the union of two cliques.
\end{obs}

For two vertices $x$ and $y$ of $G$, an $xy$-path is an induced path with ends $x$ and $y$. Throughout this paper, {\em all paths considered are induced paths}. The {\it distance} $d(x, y)$ between $x$ and $y$ is the length of the shortest $xy$-path of $G$.

Let $P$ be a path, and let $u$ and $v$ be two vertices of $P$. We use $P^*$ to denote the set of {\em internal vertices} of $P$ (i.e., those vertices of degree 2 in $P$), and use $P[u, v]$ to denote the segment of $P$ between $u$ and $v$.

Let $v\in V(G)$, and let $X$ be a subset of $V(G)$. We use $N_X(v)$ to denote the set of neighbors of $v$ in $X$. We say that $v$ is \textit{complete} to $X$ if $N_X(v)=X$, and say that $v$ is \textit{anticomplete} to $X$ if $N_X(v)=\emptyset$. For two subsets $X$ and $Y$ of $V(G)$, we say that $X$ is {\em complete} to $Y$ if each vertex of $X$ is complete to $Y$, and say that $X$ is {\em anticomplete} to $Y$ if each vertex of $X$ is anticomplete to $Y$. If $2\le|X|\le |V(G)|-1$ and every vertex in $V(G)\backslash X$ is either complete to $X$ or anticomplete to $X$, then $X$ is said to be a \textit{homogeneous set}.

\renewcommand{\baselinestretch}{1}
\begin{lemma}\label{homoset} {\em \cite{CS19}}
A minimal nonperfectly  divisible graph admits no homogeneous sets.
\end{lemma}

Let $d(v, X)=\min_{x\in X}d(v, x)$, and call $d(v, X)$  {\it the distance} of a vertex $v$ to a subset $X$. Let $i$ be a positive integer, and $N^i_G(X)=\{y\in V(G)\backslash X | d(y, X)=i\}$. We call $N^i_G(X)$ the \textit{$i$-neighborhood} of $X$, and simply write $N^1_G(X)$ as $N_G(X)$. If no confusion may occur, we write $N^i(X)$ instead of $N^i_G(X)$, and $N^i(\{v\})$ is denoted by $N^i(v)$ for short.

Suppose that  $C=v_1v_2v_3v_4v_5v_1$ is a 5-hole of $G$. For a subset $T\subseteq \{1,2,3,4,5\}$, let
$$N_T(C) = \{x\;|\; x\in N(C), \mbox{ and } v_ix\in E(G) \mbox{ if and only if } i\in T\}.$$
It is easy to check that for $k\in \{1,2,3,4,5\}$ and $l=k+2$, $N_{\{k, k+2\}}(C)=N_{\{l, l+3\}}(C)$ and $N_{\{k, k+2, k+3\}}(C)=N_{\{l, l+1, l+3\}}(C)$.

The next lemma is devoted to the structure of $P_5$-free graphs. It holds trivially by the $P_5$-freeness of $G$, and so we omit its proof.

\renewcommand{\baselinestretch}{1}
\begin{lemma}\label{P5free}
Suppose that $G$ is a $P_5$-free graph and  $C=v_1v_2v_3v_4v_5v_1$ is a $5$-hole of $G$. Then,
\begin{itemize}

\item [$(a)$] for $i\in \{1,2,3,4,5\}$, $N_{\{i\}}(C)=N_{\{i,i+1\}}(C)=\emptyset$, and $N_{\{i,i+2\}}(C)\cup N_{\{i,i+1,i+2\}}(C)$ is anticomplete to  $N^2(C)$,

\item [$(b)$] if $x\in N(C)$ and $N^2(x)\cap N^3(C)\neq\emptyset$ then $x\in N_{\{1,2,3,4,5\}}(C)$, and

\item [$(c)$] for each vertex $x\in N^2(C)$ and each component $B$ of $G[N^3(C)]$, $x$ is  either complete or anticomplete to  $B$.

\end{itemize}

\end{lemma}\renewcommand{\baselinestretch}{1.15}

We end this section by the following two lemmas which are also very useful in the proofs of the main results. A {\em clique cut set} is a cut set and is a clique.

\begin{lemma}\label{lem-clique-cut}
A minimal nonperfectly  divisible graph has no clique cut sets.
\end{lemma}\renewcommand{\baselinestretch}{1.15}
\pf If it is not the case, let $G$ be a minimal nonperfectly  divisible graph, and let $S$ be a clique cut set of $G$. Let $C_1$ be a component of $G-S$, let $G_1=G[V(C_1)\cup S]$, and let $G_2=G-V(C_1)$. Then, both $G_1$ and $G_2$ are perfectly divisible. For $i\in \{1, 2\}$, let $(A_i, B_i)$ be a perfect division of $G_i$ with $G[A_i]$ perfect and $\omega(G[B_i])<\omega(G_i)$. Since $S$ is a clique, we see that both $A_1\cap A_2$ and $B_1\cap B_2$ are cliques as they are subsets of $S$, and thus $G[A_1\cup A_2]$ is perfect and $\omega(B_1\cup B_2)<\omega(G)$, a contradiction. \qed

\medskip

Let $G$ be a graph with $\alpha(G)=2$, and let $v$ be a vertex of $G$. Notice that $V(G)\setminus(N(v)\cup\{v\})$ is a clique, which implies that $G-N(v)$ is perfect. Thus the next lemma follows directly.

\begin{lemma}\label{alpha-2}
Graphs of independent number at most $2$ are perfectly divisible.
\end{lemma}

\section{Perfect divisibility of $(P_5, C_5, K_{2,3})$-free graphs}

This section is aim to prove Theorem~\ref{main0}. A cut set $S$ is said to be a {\em minimal cut set} if any proper subset of $S$ is not a cut set of $G$. We first prove a lemma on the structure of $(P_5, C_5, K_{2, 3})$-free graphs.

\renewcommand{\baselinestretch}{1}
\begin{lemma}\label{lem-p5-c5-K23}
Suppose that $G$ is a  $(P_5, C_5, K_{2, 3})$-free graph without clique cut sets, and $S$ is a minimal cut set of $G$.  Then
\begin{itemize}
\item[$(a)$] $G-S$ has exactly two components, and for each pair of non-adjacent vertices $s_{1}, s_{2}\in S$, each $s_{1}s_{2}$-path with interior in exactly one component has length $2$,
\item[$(b)$] each vertex of $S$ is complete to at least one component of $G-S$, and
\item[$(c)$] $\alpha(G[S])=2$.
\end{itemize}
\end{lemma}\renewcommand{\baselinestretch}{1.2}
\pf Let $C_1, C_2, \ldots, C_t$ be the components of $G-S$. It is certain that $t\ge 2$. Since $S$ is a minimal cut set, we see that for each  $i\in \{1, 2, \ldots, t\}$,
\begin{equation}\label{eqa-cut}
\mbox{$N_{V(C_i)}(x)\neq \emptyset$  for each vertex $x\in S$.}
\end{equation}

Let $V_1=V(C_1)$ and $G_1=G[S\cup V_1]$,  let $G_2=G-V_1$, and let $V_2=V(G_2)\setminus S$.

Since $G$ has no clique cut set, we arbitrarily choose $s_1$ and $s_2$ to be two non-adjacent vertices  in $S$. Suppose that $G-S$ has at least 3 components, then $G_2-S$ is not connected as $G_1-S=C_1$. Let $C_2$ and $C_3$ be two components of $G_2-S$. For $i\in\{1, 2, 3\}$, let $P_i$ be an $s_1s_2$-path with interior in $C_i$ (recall that all paths considered are induced paths).

If one of $P_1, P_2$ and $P_3$ has length at least 3, then a $C_5$ or a $P_5$ appears. Otherwise, a $K_{2, 3}$ appears. Hence, $G-S$ has two components $G[V_1]$ and $G[V_2]$. This also implies that each $s_{1}s_{2}$-path with interior in $V_1$ or $V_2$ has length 2.

Let $s\in S$. It follows from (\ref{eqa-cut}) that $s$ has neighbors in both $V_1$ and $V_2$. Since both $G[V_1]$ and $G[V_2]$ are connected and $G$ is $P_5$-free,  we have that each vertex of $S$ is complete to either $V_1$  or $V_2$.

Now it is left to show that $\alpha(G[S])=2$. Suppose to its contrary that $s_3$ is a vertex in $S\setminus\{s_1, s_2\}$ anticomplete to $\{s_1, s_2\}$. Thus we have that, for each pair of $i, j\in \{1, 2\}$,  each $s_is_3$-path with interior in $V_j$ has length 2. Since $G$ induces no $K_{2, 3}$, we have that $N_{V_i}(s_1)\cap N_{V_i}(s_2)\cap N_{V_i}(s_3)=\emptyset$ for some $i\in\{1, 2\}$, and so we may assume that $N_{V_1}(s_1)\cap N_{V_1}(s_2)\cap N_{V_1}(s_3)=\emptyset$. Let $w_1\in V_1$ be a common neighbor of $s_1$ and $s_2$,  let $w_2\in V_1$ be a common neighbor of $s_2$ and $s_3$, and let $x\in V_2$ be a common neighbor of $s_1$ and $s_3$. If $w_1w_2\not\in E(G)$, then $G[\{s_{1}, w_{1}, s_{2}, w_{2}, s_{3}\}]=P_{5}$; otherwise, $G[\{s_1, s_3, w_1, w_2, x\}]=C_5$. This contradiction implies that $\alpha(G[S])=2$, which completes the proof of Lemma~\ref{lem-p5-c5-K23}. \qed

\bigskip

\noindent{\em Proof of} {\bf Theorem}~\ref{main0}. Let $G$ be a $(P_5, C_5, K_{2, 3})$-free graph. Suppose that $G$ is not perfectly divisible but every proper induced subgraph of $G$ is perfectly divisible. It is certain that $G$ is connected and not perfect. Let $S$ be a minimal cut set of $G$. By Lemma~\ref{lem-clique-cut}, $S$ is not a clique. It follows from Lemma~\ref{lem-p5-c5-K23} that $\alpha(G[S])=2$, $G-S$ has exactly two components, say $C_1$ and $C_2$, and each vertex of $S$ is either complete to $V(C_1)$ or $V(C_2)$. For $i\in \{1, 2\}$, let $V_i=V(C_i)$, and let $G_i=G[V_i\cup S]$.

Let $S_{0}\subseteq S$ be the set of vertices complete to  $V_1\cup V_2$. For $i\in \{1, 2\}$,  let $S_{i}\subseteq S\setminus S_0$ be the set of vertices only complete to $V_i$. Clearly $S=S_0\cup S_1\cup S_2$.

We claim that
\begin{equation}\label{eqa-V1-V2-clique}
\mbox{at least one of $V_1$ and $V_2$ is a clique.}
\end{equation}

Suppose to its contrary that both $V_1$ and $V_2$ are not cliques. Since $S$ is not a clique, we may choose $s_1$ and $s_2$ to be two non-adjacent vertices of $S$. Suppose that $\{s_1, s_2\}\cap S_0\neq\emptyset$. If $\{s_1, s_2\}\cap S_i\neq\emptyset$ for some $i\in\{1, 2\}$, then $V_i$ is a clique, otherwise an induced $K_{2,3}$ is obtained. Similarly, if  $\{s_1, s_2\}\subseteq S_0$, then both $V_1$ and $V_2$ must be cliques. Thus we may assume that  $\{s_1, s_2\}\cap S_0=\emptyset$.  Note that
$N_{V_i}(x)\neq \emptyset$  for each vertex $x\in S$ as $S$ is a minimal cut set.  If $\{s_1, s_2\}\subset S_1$, then $G$ induces a $K_{2, 3}$ whenever $N_{V_2}(s_1)\cap N_{V_2}(s_2)\neq\emptyset$, and $G$ induces a $P_5$ or a $C_5$ whenever $N_{V_2}(s_1)\cap N_{V_2}(s_2)=\emptyset$, both are contradictions. Similar contradiction happens if $\{s_1, s_2\}\subset S_2$. Therefore, we may suppose that $s_1\in S_1$ and $s_2\in S_2$, that is, both $S_0\cup S_1$ and $S_0\cup S_2$ are cliques. Now by Observation~\ref{obs1}, we have that $G[S]$ is perfect. Since $\omega(G-S)<\omega(G)$, it contradicts the minimal nonperfect divisibility of $G$, and which proves (\ref{eqa-V1-V2-clique}).

\medskip

Next we claim that
\begin{equation}\label{alpha-G=2}
\mbox{exact one of $V_1$ and $V_2$ is a clique}.
\end{equation}

To prove (\ref{alpha-G=2}), we will show that if $V_1$ and $V_2$ are both cliques then $\alpha(G)=2$, and hence deduce a contradiction to Lemma~\ref{alpha-2} claiming that all graphs $G$ with $\alpha(G)\le 2$ are perfectly divisible.

Suppose to its contrary that $V_1$ and $V_2$ are both cliques but $\alpha(G)>2$. Let $T=\{t_1, t_2, t_3\}$ be an independent set of $G$. It follows from Lemma~\ref{lem-p5-c5-K23} that $|T\cap S_i|=2$ and $|T\cap V_{3-i}|=1$ for some $i\in \{1, 2\}$. Without loss of generality, we assume that $t_1, t_2\in S_1$ and $t_3\in V_2$.

Note that $V_1$ and $V_2$ are both cliques, and $V_1$ is complete to $S_0\cup S_1$. If  $S_2=\emptyset$, then $(V_1\cup V_2, S_0\cup S_1)$ is a perfect division of $G$, contradicting the minimal nonperfect divisibility of $G$. Hence $S_2\neq \emptyset$.

Let $x$ be a vertex in $S_2$. Since no vertex of $S_2$ is complete to $V_1$, we may choose a vertex, say $v_1$, in $V_1$ with $xv_1\not\in E(G)$. Since $\alpha(G[S])=2$, we have that $x$ cannot be anticomplete to $\{t_1, t_2\}$. Suppose $xt_1\in E(G)$. To avoid a $P_5=t_2v_1t_1xt_3$, we have that $x$ must be adjacent to $t_2$ as well. Hence we have that
$\{t_1, t_2\}$ is complete to $S_2$.

If $S_2$ is not a clique, let $x$ and $x'$ be two non-adjacent vertices of $S_2$, then $G[\{T\cup \{x, x'\}]$ is a $K_{2, 3}$. This implies that $S_2$ must be a clique.

Since both $V_1$ and $S_2\cup V_2$ are cliques and $V_1$ is complete to $S_0\cup S_1$,  we have that $G[V_1\cup V_2\cup S_2]$ is perfect by Observation~\ref{obs1}, and $\omega(G[S_0\cup S_1])<\omega(G)$. Thus $(V_1\cup V_2\cup S_2, S_0\cup S_1)$ is a perfect division of $G$, which leads to a contradiction  and proves (\ref{alpha-G=2}).

\medskip

Now we may assume that $\mbox{$V_1$ is a clique and $V_2$ is not.}$

Since $V_2$ is not a clique, we must have that $S_0\cup S_2$ is a clique, otherwise an induced $K_{2,3}$ appears. Thus  $V_1\cup S_0$ is also a clique.  Since $G$ is $(P_5, C_5, K_{2, 3})$-free,  by Observation~\ref{obs1} we have that
\begin{equation}\label{eqa-V1-S0-S2}
\mbox{$G[V_1\cup S_0\cup S_2]$ is perfect}.
\end{equation}

Suppose that $S_0\cup S_2\neq \emptyset$, and let $v\in S_0\cup S_2$. If $\omega(G[V_2\cup S_1])< \omega(G)$, then $(V_1\cup S_0\cup S_2, V_2\cup S_1)$ is a perfect division of $G$ by (\ref{eqa-V1-S0-S2}). Thus we may further assume that $\omega(G[V_2\cup S_1])= \omega(G)$. Let $K_1, \ldots, K_r$ be all the cliques of $G[V_2\cup S_1]$ of size  $\omega(G)$, and let $L_i=K_i\cap S_1$ for each $i\in \{1, 2, \ldots, r\}$. Since $|K_i|=\omega(G)$ and $\omega(G[V_2])<\omega(G)$, we have that $L_i\neq \emptyset$, and $v$ is not complete to $L_i$. Let $M_i\subseteq L_i$ be the set of vertices which are non-adjacent to $v$ for $i\in \{1, 2, \ldots, r\}$, and let $M=\bigcup\limits_{i=1}^{r} M_i$. Since $\alpha(G[S])=2$, we have that $M$ is a clique, and so  $V_1\cup M$ is a clique. Notice that $S_0\cup S_2$ is a clique. Thus  $G[V_1\cup S_0\cup S_2\cup M]$ induces no odd antihole by Observation~\ref{obs1}, and so it is perfect. Now we have that $G$ is perfectly divisible as $\omega(G[(V_2\cup S_1)\setminus M])<\omega(G[V_2\cup S_1])=\omega(G)$, which contradicts the minimal nonperfect divisibility of $G$.

\medskip

Hence we may suppose that, for each minimal cut set $S$ of $G$,
$S_0\cup S_2=\emptyset$. Consequently, we must have that $|V_1|=1$ (as otherwise $V_1$ is a homogeneous set, contradicting Lemma~\ref{homoset}), that is,
\begin{equation}\label{eqa-S0-S2-emptyset}
\mbox{every minimal cut set of $G$ equal $N(x)$ for some vertex $x$ of $G$.}
\end{equation}

Let $v$ be a vertex of $G$. It is certain that $N(v)$ is a cut set. Next we show that
\begin{equation}\label{eqa-neighbors}
\mbox{$N(v)$ is a minimal cut set,}
\end{equation}
which implies that the converse of (\ref{eqa-S0-S2-emptyset}) holds as well.

Suppose to its contrary that $w$ is a vertex such that $N(w)$ is not a minimal cut set. Let $T_1, T_2, \ldots, T_r$ be all the subsets of $N(w)$ where each one is a minimal cut set. It follows from (\ref{eqa-S0-S2-emptyset}) that there are some vertices, say $w_1, w_2, \ldots, w_r$, such that $T_i=N(w_i)$ for each $i\in\{1, 2, \ldots, r\}$. We claim that
\begin{equation}\label{eqa-all-single-vertex}
\{w, w_1, w_2, \ldots, w_r\}=V(G)\setminus N(w).
\end{equation}

If it is not the case, then let $C$ be a component of $G- N(w) - \{w, w_1, w_2, \ldots, w_r\}$, and let $X\subseteq N(w)$ be the set of vertices where each one has a neighbor in $C$. It is certain that $X$ is a cut set. Thus there must be some $i$ such that $T_i\subseteq X$. Without loss of generality, we suppose that $T_1\subseteq X$. Since $G$ has no clique cut set by Lemma~\ref{lem-clique-cut}, we have that $T_1$ is not a clique, and so has two non-adjacent vertices, say $t_1$ and $t'_1$. Let $P$ be a shortest $t_1t'_1$-path with interior in $C$. If $P$ has length 2, then $t_1wt'_1, t_1w_1t'_1$ and $P$ form an induced $K_{2, 3}$ in $G$. If $P$ has length 3, then $t_1wt'_1$ and $P$ form a $C_5$ in $G$. Otherwise, we have that $P$ has length greater than 3, and then we can find a $P_5$ in $G$. These contradictions proves (\ref{eqa-all-single-vertex}).

Since $\{w, w_1, w_2, \ldots, w_r\}$ is independent, it follows from (\ref{eqa-all-single-vertex}) that $(\{w, w_1, w_2, \ldots, w_r\}, N(w))$ forms a perfect division of $G$, and so (\ref{eqa-neighbors}) holds. Thus by Lemma~\ref{lem-p5-c5-K23}, we have that
\begin{equation}\label{eqa-neighbors-single}
\mbox{$\alpha(G[N(x)])=2$ for each vertex $x$ of $G$.}
\end{equation}

We choose $v_1\in V(G)$ and let $S_1=N(v_1)$. Then, $\alpha(G[S_1])=2$, and $S_1$ is a minimal cut set by (\ref{eqa-neighbors}). Let $s_1$ and $s_2$ be two non-adjacent vertices in $S_1$.  Let $V_1=\{v_1\}$,  and let $V_2=V(G)\setminus \{S_1\cup \{v_1\})$. We have that $G[V_2]$ is connected by Lemma~\ref{lem-p5-c5-K23}$(a)$.

Let $M=N_{V_2}(s_1)\cap N_{V_2}(s_2)$, and let $M_i=N_{V_2}(s_i)\setminus M$ for $i\in \{1, 2\}$. Since $G$ induces no $K_{2, 3}$, we have that $M$ must be a clique. If both $M_1$ and $M_2$ are not empty, let $m_i\in M_i$ for $i\in \{1, 2\}$, then $G[\{m_1, m_2, s_1, s_2, v_1\}]=C_5$ whenever $m_1m_2\in E(G)$, and $G[\{m_1, m_2, s_1, s_2, v_1\}]=P_5$ whenever $m_1m_2\not\in E(G)$. Without loss of generality, we assume that $M_2=\emptyset$. Now we claim that
\begin{equation}\label{Mv1=M1+M}
V_2=M\cup M_1.
\end{equation}

Suppose that (\ref{Mv1=M1+M}) does not hold. Since $G[V_2]$ is connected, we may choose a vertex, say $z$, in $V_2\setminus(M\cup M_1)$ that is adjacent to some vertex of $M\cup M_1$. If $zz_1\in E(G)$ for some $z_1\in M$, then $\{s_1, s_2, z\}$ is an independent set contained in $N(z_1)$, contradicting (\ref{eqa-neighbors-single}). If $zz_2\in E(G)$ for some vertex $z_2\in M_1$, then $G[\{s_1, s_2, v_1, z, z_2\}]=P_5$. This proves (\ref{Mv1=M1+M}).

\medskip

Note that $S_1=N(v_1)$. By (\ref{eqa-neighbors-single}), we have that no vertex of
$N(v_1)\setminus \{s_1, s_2\}$ is anticomplete to $\{s_1, s_2\}$. Let $T=N_{S_1}(s_1)\cap N_{S_1}(s_2)$, let $T_1=N_{S_1}(s_1)\setminus T$,  and let $T_2=N_{S_1}(s_2)\setminus T$. If $T_i$ is not a clique for some $i\in \{1, 2\}$, let $t_{i, 1}$ and $t_{i, 2}$ be two non-adjacent vertices of $T_i$, then $\{s_{3-i}, t_{i, 1}, t_{i, 2}\}$ is an independent set in $S_1$, which leads to a contradiction to (\ref{eqa-neighbors-single}). Thus both $T_1$ and $T_2$ are cliques.

Let $A=T_2\cup\{s_1, s_2\}$, and let $B=V(G)\setminus A$. It is certain that $G[A]$ is perfect. Since $s_1$ is complete to $B$ by (\ref{Mv1=M1+M}), we see that $\omega(G[B])<\omega(G)$, which implies that $(A, B)$ is a perfect division of $G$. This contradicts the minimal nonperfect  divisibility of $G$ and proves Theorem~\ref{main0}. \qed

\section{$(P_5, K_{2,3})$-free graphs}

In this section, we prove Theorem~\ref{main1}.

Since perfectly divisible graphs $G$ has chromatic number at most ${\omega(G)+1\choose 2}$, it follows from Theorem~\ref{main0} that we only need to consider the $(P_5, K_{2,3})$-free graphs with a 5-hole. Let $G$ be a $(P_5, K_{2,3})$-free graph, and let $C=v_1v_2v_3v_4v_5v_1$ be a 5-hole of $G$. Recall that for $T\subseteq \{1,2,3,4,5\}$,
$N_T(C)$ consists of the vertices not on $C$ but each has exactly $\{v_i\;|\; i\in T\}$ as its neighbors on $C$, and for integer $i\ge 1$, $N^i(C)$ consists of the vertices of distance $i$ apart from $C$.
Let $u$ and $v$ be two non-adjacent vertices in $N(C)$. We say that $\{u, v\}$ is a {\em bad pair} if there is an $i\in \{1,2,3,4,5\}$ such that $u\in N_{\{i,i+1,i+3\}}(C)$ and $v\in N_{\{i,i+1,i+2,i+4\}}(C)$.

\renewcommand{\baselinestretch}{1}
\begin{lemma}\label{P5K23free}
Suppose that $G$ is a $(P_5, K_{2,3})$-free graph with a $5$-hole $C=v_1v_2v_3v_4v_5v_1$, and $u, v$ are two vertices in $N(C)$. Then all the followings hold.
\begin{itemize}

\item [$(a)$] If there exist three consecutive vertices of $C$, named $v_i, v_{i+1}$, and $v_{i+2}$, such that $\{v_i, v_{i+2}\}\subseteq N(u)\cap N(v)$ and $v_{i+1}\notin N(u)\cup N(v)$, then $uv\in E(G)$.

\item [$(b)$] $N_{\{i, i+2\}}(C)$, $N_{\{i, i+1, i+3\}}(C)$, and $N_{\{i,i+1,i+2,i+3\}}(C)$ are all cliques  for $1\le i\le 5$.

\item [$(c)$] $\alpha(G[N_{\{1,2,3,4,5\}}(C)])\le 2$, and for each $i\in \{1, 2, 3, 4, 5\}$, $\alpha(G[N_{\{i,i+1,i+2\}}(C)])\le 2$, and $N_{\{i, i+1, i+2\}}(C)$ is complete to $N_{\{i+1, i+2, i+3\}}(C)$.

\item [$(d)$] If $uv\not\in E(G)$ and $\{u, v\}$ is not a bad pair, then $N(u)\cap N(v)\cap N^2(C)=\emptyset$.

\item [$(e)$] $N(C)\backslash (N_{\{1,2,3,4,5\}}(C)\bigcup\limits_{1\le i\le5} N_{\{i,i+1,i+2\}}(C))$ can be partitioned into five cliques $S_1, S_2, S_3, S_4$ and $S_5$ such that $|S_i|\le \omega(G)-1$ and $v_i$ is anticomplete to $S_i$ for each $i$.

\end{itemize}
\end{lemma}\renewcommand{\baselinestretch}{1.2}
\pf If $\{v_i, v_{i+2}\}\subseteq N(u)\cap N(v)$ and $v_{i+1}\notin N(u)\cup N(v)$ for some $i$, then $uv\in E(G)$ to avoid a $K_{2,3}$ on $\{u, v, v_{i}, v_{i+1}, v_{i+2}\}$. Hence $(a)$ holds, and $(b)$ follows directly from $(a)$.

Now we come to prove $(c)$. If either $G[N_{\{i,i+1,i+2\}}(C)]$ or $G[N_{\{1,2,3,4,5\}}(C)]$  has an independent set of size 3, say $\{u, v, w\}$, then $\{u, v, w, v_{i}, v_{i+2}\}$ induces a $K_{2,3}$. If $N_{\{i,i+1,i+2\}}(C)$ is not complete to $N_{\{i+1,i+2,i+3\}}(C)$ for some $i$, we may choose $x\in N_{\{i,i+1,i+2\}}(C)$ and $y\in N_{\{i+1,i+2,i+3\}}(C)$ with $xy\not\in E(G)$,  then a $P_5=xv_{i+1}yv_{i+3}v_{i+4}$ appears. Hence $(c)$ holds.

\medskip

Next we prove $(d)$. Suppose that $uv\not\in E(G)$ and  $\{u, v\}$ is not a bad pair, and suppose that $N(u)\cap N(v)\cap N^2(C)$ has a vertex $w$. If one of $u$ and $v$ is in $N_{\{1,2,3,4,5\}}(C)$, then there exists some $i\in \{1,2,3,4,5\}$ such that $G[\{u, v, v_i, v_{i+2}, w\}]=K_{2, 3}$, which leads to a contradiction.  Thus $\{u, v\}\subseteq W=\bigcup\limits_{1\le i\le 5}(N_{\{i, i+1, i+3\}}(C)\cup N_{\{i,i+1,i+2,i+3\}}(C))$ by Lemma~\ref{P5free}$(a)$.

Now we first suppose that $u\in N_{\{k, k+1, k+3\}}(C)$ for some $k$, and by symmetry we may assume that $k=1$.
Since $\{u, v\}$ is not a bad pair, we see that $v\in W\setminus N_{\{1,2,3,5\}}(C)$. It follows from $(a)$ that, $v\notin N_{\{1, 2, 4\}}(C)\cup N_{\{1, 3, 4\}}(C)\cup N_{\{2, 4, 5\}}(C)\cup N_{\{1,2,3,4\}}(C)\cup N_{\{4,5,1,2\}}(C)$. But $G[\{u, v, w, v_2, v_5\}]=P_5$ if $v\in N_{\{1,3,5\}}(C)$, $G[\{u, v, w, v_1, v_3\}] = P_5$ if $v\in N_{\{2,3,5\}}(C)\cup N_{\{2,3,4,5\}}(C)$, and $G[\{u, v, w, v_1, v_4\}] = P_5$ if $v\in N_{\{3,4,5,1\}}(C)$. Hence we have, by symmetry, that
$$\{u, v\}\cap (\bigcup\limits_{1\le i\le 5}(N_{\{i, i+1, i+3\}}(C))=\emptyset,$$
that is, $\{u,v\}\subseteq \bigcup\limits_{1\le i\le 5} N_{\{i,i+1,i+2,i+3\}}(C)$.
Without loss of generality, we may assume that $u\in N_{\{1,2,3,4\}}(C)$. Thus $v\notin \bigcup\limits_{1\le i\le 5} N_{\{i,i+1,i+2,i+3\}}(C)$ by $(a)$ of this lemma. This contradiction proves $(d)$.

Finally, we prove $(e)$. Let $S_1=N_{2,5}(C)\cup N_{\{2,3,5\}}(C)\cup N_{\{2,4,5\}}(C)\cup N_{\{2,3,4,5\}}(C)$,
$S_2=N_{\{1,3\}}(C)\cup N_{\{1,3,4\}}\cup N_{\{1,3,5\}}\cup N_{\{1,3,4,5\}}$, $S_3=N_{\{2,4\}}\cup N_{\{1,2,4,5\}}(C)$, $S_4=N_{\{3, 5\}}(C) \cup  N_{\{1,2,3,5\}}(C)$, and $S_5=N_{\{1,4\}}(C)\cup N_{\{1,2,4\}}(C)\cup N_{\{1,2,3,4\}}(C)$. It is certain that $v_i$ is anticomplete to $S_i$ for each $i$, and one can check easily from $(a)$ that  each $S_i$ is a clique of size at most $\omega(G)-1$.  \qed

\renewcommand{\baselinestretch}{1}
\begin{lemma}\label{P5K23freecliguecut set}
Let $G$ be a $(P_5, K_{2,3})$-free graph and $C=v_1v_2v_3v_4v_5v_1$ be a $5$-hole of $G$. If $G$ contains no clique cut sets, then $N^3(C)=\emptyset$, and for each component $B$ of $N^2(C)$, $\alpha(B)\le 2$ and $N(B)\cap N(C)\subseteq N_{\{1, 2, 3, 4, 5\}}(C)$ whenever $\omega(B)=\omega(G)$.
\end{lemma}\renewcommand{\baselinestretch}{1.2}
\pf Suppose that $G$ has no clique cut sets. Let $S'\subseteq N(C)$ be the set of vertices  having no neighbors in $N^2(C)$, and let $S=N(C)\backslash S'$. By Lemma~\ref{P5free}$(a)$, if a vertex $x$ of $N(C)$ has neighbors in $N^2(C)$, then $|N_{C}(x)|\ge 3$.  Hence for any two vertices $u$ and $v$ of $S$ we have that
\begin{equation}\label{eqa-0}
N_{C}(v)\cap N_{C}(u)\neq\emptyset.
\end{equation}

Next we prove this lemma by considering the connectedness of $G[N^2(C)]$.

\medskip

\noindent{\bf Case 1}. Suppose that $G[N^2(C)]$ is connected. Since $G$ contains no clique cut sets, there must exist two non-adjacent vertices, say  $u$ and $v$, in $S$. By (\ref{eqa-0}), we may choose $z$ to be a common neighbor of $u$ and $v$ on $C$.

Let $T_{uv}=N(u)\cap N(v)\cap  N^2(C)$, and let $T_u=(N(u)\cap N^2(C))\setminus T_{uv}$ and $T_v=(N(v)\cap  N^2(C))\setminus T_{uv}$.

\medskip

First suppose that $T_{uv}=\emptyset$. Since $G$ is $(P_5, K_{2, 3})$-free, we have
\begin{equation}\label{eqa-1}
\mbox{$T_u$ is complete to $T_v$, and $\max\{\alpha(G[T_u]), \alpha(G[T_v])\le 2$},
\end{equation}
otherwise for any pair of non-adjacent vertices $t_u\in T_u$ and $t_v\in T_v$, $t_uuzvt_v$ is a $P_5$, and $G[T_u\cup T_v\cup \{w\}]$ induces a $K_{2, 3}$ whenever $\alpha(G[T_w])\ge 3$ for any $w\in \{u, v\}$. This leads to a contradiction.

We further claim that

\begin{equation}
N^2(C)=T_u\cup T_v, \mbox{ and } N^3(C)=\emptyset. \label{eqa-2}
\end{equation}

If $N^2(C) \neq T_u\cup T_v$, since $G[N^2(C)]$ is connected, we may suppose by symmetry that $N^2(C)\backslash(T_u\cup T_v)$ has a vertex, say $w$, which has a neighbor $t_u$ in $T_u$, and so $vzut_uw$ is a $P_5$. Thus $N^2(C) = T_u\cup T_v$. If $N^3(C)\neq \emptyset$, then let $w$ be a vertex of $N^3(C)$, and suppose by symmetry that $t'_u$ is a neighbor of $w$ in $T_u$, again we have a $P_5=vzut'_uw$. This proves (\ref{eqa-2}).

From (\ref{eqa-1}) and (\ref{eqa-2}), we see that $\alpha(G[N^2(C)])\le 2$. If $\omega(G[N^2(C)])=\omega(G)$, then $T_u\neq\emptyset$ and $T_v\neq \emptyset$. To avoid a $P_5$, we see that all neighbors of $N^2(C)$ in $N(C)$ must be contained in $N_{\{1, 2, 3, 4, 5\}}(C)$.
So the lemma holds whenever $G[N^2(C)]$ is connected and $T_{uv}=\emptyset$.

\medskip

Next we may assume that $T_{uv}\neq\emptyset$, and let $t_{uv}$ be a vertex in $T_{uv}$. By Lemma~\ref{P5K23free}$(d)$, we have that $\{u, v\}$ is a bad pair. Suppose by symmetry that $u\in N_{\{1,2,4\}}(C)$ and $v\in N_{\{1, 2, 3, 5\}}(C)$.  Since $u$ and $v$ have a common neighbor $z$ on $C$, and since $G$ is $K_{2,3}$-free, we have that
\begin{equation}\label{eqa-Tuv-clique}
\mbox{$T_{uv}$ is a clique.}
\end{equation}

Suppose that $T_v\neq\emptyset$. Let $t_v$ be a vertex of $T_v$. If $t_v$ has a neighbor $t'_{uv}$ in $T_{uv}$, then a $P_5=v_{5}v_1ut'_{uv}t_v$ appears, and if $t_v$ is anticomplete to $T_{uv}$, then a $P_5=v_{4}ut_{uv}vt_v$ appears. This implies that $T_v=\emptyset$. Similarly we have that $T_u=\emptyset$.

Note that $G[N^2(C)]$ is connected. Thus $N^2(C)=T_{uv}$, otherwise we may choose a vertex $w$ in $N^2(C)\setminus T_{uv}$ where $w$ has a neighbor $t_{uv}\in T_{uv}$ which implies a $P_5=v_{5}v_1ut_{uv}w$. If $N^3(C)\ne \emptyset$, let $w$ be a vertex in $N^3(C)$ and $w'$ be a neighbor of $w$ in $N^2(C)$, then a $P_5=ww'uv_1v_5$ appears. Thus $N^3(C)=\emptyset$, and by (\ref{eqa-Tuv-clique}) we have that this lemma holds when $G[N^2(C)]$ is connected and $T_{uv}\neq\emptyset$. Hence Lemma~\ref{P5K23freecliguecut set} holds when $G[N^2(C)]$ is connected.

\medskip

\noindent{\bf Case 2}. Suppose that  $G[N^2(C)]$ is not connected. Let $T_1, T_2, \ldots, T_k$ be the components of $G[N^2(C)]$. Recall that $S$ consists of the vertices  in $N(C)$ that have neighbors in $N^2(C)$.   For each $i\in \{1,2,\ldots, k\}$, let $S_i=N_S(T_i)$, and let $Z_i$ be the subgraph of $G[N^3(C)]$ that is not anticomplete to $T_i$. If there exists some $i_0\in \{1,2,\ldots, k\}$ such that $S_{i_0}$ is not a clique, by applying the same arguments to $G[V(C)\cup S_{i_0}\cup T_{i_0}\cup Z_{i_0}]$ as that used in Case 1, we can show that $Z_{i_0}=\emptyset$ and $\alpha(T_{i_0})\le 2$, and $N(T_{i_0})\cap N(C)\subseteq N_{\{1, 2, 3, 4, 5\}}(C)$ whenever $\omega(T_{i_0})=\omega(G)$. Hence the lemma holds if $S_i$ is not a clique for all $1\le i\le k$.

Thus by symmetry suppose that $S_1$ is a clique. Since $S_1$ is not a  clique cut set, we have that $S_1\neq S$. Without loss of generality, let $T_1, T_2, \cdots, T_{\l}$ be the components such that $N_S(T_i)\subseteq S_1$ for each $i\in \{1, 2, \ldots, \l\}$. It is obvious that $1\le \l\le k-1$, otherwise $S_1$ is a clique cut set. Now we have that $T_i$ has neighbors in $S\backslash S_1$ for each $i\in \{\l+1, \ldots, k\}$.

Since $S_1$ is not a clique cut set, we have that $N^3(C)\neq\emptyset$, and $T_1$ must have some neighbors in $N^3(C)$. Let $R$ be a component of $G[N^3(C)]$ such that $T_1$ is not anticomplete to $R$. Since $S_1$ is not a clique cut set, we have that $R$ cannot be anticomplete to $\cup_{i=l+1}^k T_i$. Without loss of generality, suppose that $R$ is not anticomplete to $T_k$. Choose $t_1\in T_1$, $t_k\in T_k$, and  $r_1, r_k\in R$ such that $t_1r_1\in E(G)$ and $t_kr_k\in E(G)$. By Lemma~\ref{P5free}$(c)$, $\{t_1, t_k\}$ is complete to $R$.

We can choose two adjacent vertices $s_k\in S\backslash S_1$ and $t'_k\in T_k$. Let $P'$ be a shortest $t_kt'_k$-path in $T_k$. Let $r\in R$ and $z$ be a neighbor of $s_k$ on $C$. Thus a path $P=t_1rt_kP't'_ks_kz$ of length at least 5 appears. This contradiction proves Lemma~\ref{P5K23freecliguecut set}.  \qed

\medskip

Now, we can prove Theorem~\ref{main1}.

\medskip

\noindent\textbf{Proof of Theorem~\ref{main1}.} Let $G$ be a $\{P_5, K_{2,3}\}$-free graph. We may assume that $G$ is connected and contains no clique cut set. If $G$ is $(P_5, C_5, K_{2,3})$-free, then $G$ is perfectly divisible by Theorem~\ref{main0}, which implies that $\chi(G)\le \frac{1}{2}(\omega^2(G)+\omega(G))$. Thus we suppose that $G$ is $(P_5, K_{2, 3})$-free and contains a $5$-hole $C$. By Lemma~\ref{P5free}$(a)$, we have that
\begin{equation*}
\begin{aligned}
N(C)  & =N_{\{1,2,3,4,5\}}(C)\bigcup_{1\le i\le 5} (N_{\{i, i+2\}}(C)\cup N_{\{i, i+1, i+2\}}(C) \\
          & \ \ \ \cup N_{\{i, i+1, i+3\}}(C)\cup N_{\{i, i+1, i+2, i+3\}}(C)).
\end{aligned}
\end{equation*}
Note that $\omega(G[N_{\{1,2,3,4,5\}}(C)])\le \omega(G)-2$, $\omega(G[N_{\{1, 2, 3\}}(C)\cup N_{\{2, 3, 4\}}(C)]\le \omega(G)-2$,  $\omega(G[N_{\{3, 4, 5\}}(C)\cup N_{\{1, 4, 5\}}(C)]\le \omega(G)-2$, and $\omega(G[N_{\{1, 2, 5\}}(C)]\le \omega(G)-2$.
By Lemma~\ref{alpha-2} and Lemma~\ref{P5K23free}$(c)(e)$,
\begin{eqnarray}\label{eqa-NC-0}\nonumber
\chi(G[N(C)]) & \le & \chi(G[N_{\{1, 2, 3\}}(C)\cup N_{\{2, 3, 4\}}(C)]) + \chi(G[N_{\{3, 4, 5\}}(C)\cup N_{\{1, 4, 5\}}(C)])\\\nonumber
     &  & + \chi(G[N_{\{1,2,5\}}(C)])+ \chi(G[N_{\{1,2,3,4,5\}}(C)])+ 5(\omega(G)-1)\\\nonumber
     & \le & 4\cdot\frac{(\omega(G)-2)^2+\omega(G)-2}{2}+5(\omega(G)-1)\\
     & = & 2\omega^2(G)-\omega(G)-3.
\end{eqnarray}

By Lemma~\ref{P5K23free}$(e)$, we can color the vertices of $C$ with the colors used on the vertices of $N(C)\backslash (N_{\{1,2,3,4,5\}}(C)\bigcup\limits_{1\le i\le5} N_{\{i,i+1,i+2\}}(C))$ (which is counted  $5(\omega(G)-1)$ in (\ref{eqa-NC-0})).

Let $B$ be a component of $G[N^2(C)]$. By Lemma~\ref{P5K23freecliguecut set}, we see that $N^2(C)=G-N(C)-V(C)$, $\alpha(B)\le 2$ and $N(B)\cap N(C)\subseteq N_{\{1, 2, 3, 4, 5\}}(C)$ if $\omega(B)=\omega(G)$. So, we have, by Lemmas~\ref{alpha-2}, that $\chi(B)\le \frac{(\omega(G)-1)^2+\omega(G)-1}{2}=\frac{\omega^2(G)-\omega(G)}{2}$ if $\omega(B)<\omega(G)$, and $\chi(B)\le \frac{\omega^2(G)+\omega(G)}{2}$ otherwise.

Note that $N^2(C)$ is anticomplete to  $\bigcup\limits_{1\le i\le 5} N_{\{i, i+1, i+2\}}(C)$ by Lemma~\ref{P5free}$(a)$, and $B$ is anticomplete to $N(C)\backslash N_{\{1,2,3,4,5\}}(C)$  by Lemma~\ref{P5K23freecliguecut set} if $\omega(B)=\omega(G)$. If $\omega(B)<\omega(G)$, we can color the vertices in $B$ with the colors used on the vertices of $\bigcup\limits_{1\le i\le 5} N_{\{i, i+1, i+2\}}(C)$ (which is counted no less than $\frac{\omega^2(G)-\omega(G)}{2}$ in (\ref{eqa-NC-0})). If $\omega(B)=\omega(G)$, we can color the vertices in $B$ with the colors used on the vertices of $N(C)\backslash N_{\{1,2,3,4,5\}}(C)$ (which is counted no less than $\frac{\omega^2(G)+\omega(G)}{2}$ in (\ref{eqa-NC-0})).

Therefore, $\chi(G)\le 2\omega^2(G)-\omega(G)-3$ as desired. \qed

\section{$(P_5, K_1+2K_2)$-free graphs}

For two subsets $X$ and $Y$ of $V(G)$, we say that $X$ {\em dominates} $Y$ if each vertex of $Y$ has a neighbor in $X$. The next two lemmas are very useful in the proof of Theorem~\ref{main2}.

\renewcommand{\baselinestretch}{1}
\begin{lemma}{\em \cite{BT90}}\label{P5dom}
Every connected $P_5$-free graph has a dominating clique or a dominating $P_3$.
\end{lemma}

\begin{lemma}{\em \cite{W80}}\label{2K2free}
Let $G$ be a $2K_2$-free graph. Then  $\chi(G)\le \frac{1}{2}(\omega^2(G)+\omega(G))$.
\end{lemma} \renewcommand{\baselinestretch}{1.2}

\noindent\textbf{Proof of Theorem~\ref{main2}.} Let $G$ be a connected $(P_5, K_1+2K_2)$-free graph with at least two vertices. By Lemma~\ref{P5dom}, $G$ has a  dominating clique or a dominating $P_3$. If $G$ has a dominating $P_3$, say $v_1v_2v_3$, then $N(v_i)$ induces a $2K_2$-free graph for each $i$, otherwise $v_i$ and the $2K_2$ in $G[N(v_i)]$ induce a $K_1+2K_2$. Thus by Lemma~\ref{2K2free} we have that $\chi(G)\le \chi(G[N(v_1)])+\chi(G[N(v_2)])+\chi(G[N(v_3)])\le \frac{3}{2}((\omega(G)-1)^2+\omega(G)-1)={3\over 2}(\omega^2(G)-\omega(G))$. Thus we may assume that $G$ has a dominating clique, say $K_k$ on vertices $\{v_1, v_2, v_3, \cdots, v_k\}$. Let $S=N(v_1)\cup N(v_2)$ and $T=V(G)\backslash S$, that is, $T$ consists of exactly those vertices anticomplete to $\{v_1, v_2\}$. Since $G$ is $(K_1+2K_2)$-free, we have that $N(v_i)$ induces a $2K_2$-free subgraph for $i=1,2$. For each $i\in \{3, \dots, k\}$, since the vertices of $N(v_i)\cap T$ are anticomplete to $\{v_1, v_2\}$, we have that $N(v_i)\cap T$ is an independent set. Hence by Lemma~\ref{2K2free}, $\chi(G)\le \chi(G[N(v_1)])+\chi(G[N(v_2)])+\omega(G)-2\le \omega^2(G)-2$. Note that ${3\over 2}(\omega^2(G)-\omega(G))\ge \omega^2(G)-2$. Hence if $G$ is $(P_5, K_1+2k_2)$-free, then $\chi(G)\le {3\over 2}(\omega^2(G)-\omega(G))$ as desired. \qed

\section{$(P_5, K_1+(K_1\cup K_3))$-free graphs}

Sumner \cite{sumner1} (see also \cite{ELMM13}) proved that a connected $(P_5, K_3)$-free graph is either bipartite or can be obtained from a 5-hole by replacing each vertex with an independent set and then replacing each edge by a complete bipartite graph.
\begin{lemma}\label{P5K3} {\em \cite{sumner1}}
If $G$ is  $(P_5, K_3)$-free then  $\chi(G)\le 3$.
\end{lemma}

If $G$ is $(P_5, K_1\cup K_3)$-free, then $G-N(v)-v$ is $(P_5, K_3)$-free for each vertex $v$ of $G$, and so by a simple induction one can show the following lemma.

\begin{lemma} \label{P5K1K3}
$\chi(G)\le 3\omega(G)-3$ for every $(P_5, K_1\cup K_3)$-free graph $G$ with at least one edge.
\end{lemma}

Before proving Theorem~\ref{main3}, we first prove a few lemmas on the structure of $(P_5, K_1+(K_1\cup K_3))$-free graphs. From now on, we always suppose that $G$ is a $(P_5, K_1+(K_1\cup K_3))$-free graph.

\renewcommand{\baselinestretch}{1}
\begin{lemma}\label{P5K1K3K1free}
Suppose that  $G$ has a $5$-hole $C=v_1v_2v_3v_4v_5v_1$  and has no clique cut set, and let $T$ be a component of $G[N^2(C)]$. Then the followings hold.
\begin{itemize}

\item [$(a)$] For each $i\in \{1,2,3,4,5\}$, $G[N(v_i)]$ is $K_1\cup K_3$-free, $G[N_{\{i,i+2\}}(C)]$ is $K_3$-free,  and $N_{\{i,i+1,i+2\}}(C)\cup N_{\{i,i+1,i+3\}}(C)\cup N_{\{i,i+1,i+2,i+3\}}(C)$ is independent.

\item [$(b)$] If no vertex in $N(C)$ dominates $T$, then there exist two non-adjacent vertices $u$ and $v$ in $N(C)$ such that both $N_T(u)$ and $N_T(v)$ are not empty.

\end{itemize}
\end{lemma}\renewcommand{\baselinestretch}{1.2}
\pf Statement $(a)$ follows directly from the $K_1+(K_1\cup K_3)$-freeness of $G$.

To prove $(b)$, let $S= N(T)\cap N(C)$, and suppose that no vertex in $S$ dominates $T$. By Lemma~\ref{P5free}$(a)$, $S\subseteq N_{\{1,2,3,4,5\}}(C)\cup (\bigcup\limits_{1\le i\le 5}N_{\{i,i+2,i+3\}}(C)\cup N_{\{i,i+1,i+2,i+3\}}(C))$. If $G[S]$ is not connected, we are done. Thus suppose that $G[S]$ is connected.  We choose an arbitrary vertex $u\in S$ and let $T_u=N_T(u)$. Since $u$ does not dominate $T$, by the connectedness of $T$, we may choose a vertex $w\in T\backslash T_u$ such that $w$ is not anticomplete to $T_u$. Let $v$ be a neighbor of $w$ in $S$. Since $G$ is $P_5$-free, we have that $u\in N_{\{1,2,3,4,5\}}(C)$, and so $\{v_{i+2}, v_{i+3}\}\subseteq N(u)\cap N(v)$. Since $uv\in E(G)$ implies a $K_1+(K_1\cup K_3)$ on  $\{w,u,v,v_{i+2},v_{i+3}\}$, we have that $uv\notin E(G)$ as desired.   \qed

\renewcommand{\baselinestretch}{1}
\begin{lemma}\label{P5K1K3K1freeC5}
Suppose that $G$ has a $5$-hole $C=v_1v_2v_3v_4v_5v_1$ and no clique cut set. Then $G[N^3(C)]$ is $K_3$-free, and  $N^2(C)$ can be partition into two parts $A$ and $B$  such that both  $G[A]$ and $G[B]$ are $K_3$-free.
\end{lemma}\renewcommand{\baselinestretch}{1.2}
\pf Let $B$ be a component of $G[N^3(C)]$ and $u\in N^2(C)$ be a vertex that has a neighbor in $B$. By Lemma~\ref{P5free}$(c)$, we see that $u$ must be complete to $B$, and so $G[N^3(C)]$ must be $K_3$-free  to avoid a $K_1+(K_1\cup K_3)$.

Let $T=N^2(C)$. Without loss of generality, we suppose that $G[T]$ is connected, and let
$$S=\{v|v\in N(C) \mbox{ such that } N_T(v)\neq\emptyset\}.$$

If  there exists  some vertex in $S$   that dominates $T$, then we are done as $G[T]$ is obviously $K_3$-free to avoid a $K_1+(K_1\cup K_3)$. Thus suppose that no vertex of $S$ dominates $T$.

By Lemma~\ref{P5K1K3K1free}$(b)$, there exist two non-adjacent vertices, say  $u$ and $v$, in $S$ such that both $u$ and $v$ have neighbors  $T$.  It follows from Lemma~\ref{P5free}$(a)$ that $u$ and $v$ have a common neighbor, say $z$, on $C$.

It is certain that both $G[N_T(u)]$ and $G[N_T(v)]$ are $K_3$-free. If $T=N(u)\cup N(v)$, then $(N_T(u), N_T(v)\setminus N_T(u))$ is a partition of $T$ as desired. Thus suppose that $T\neq N(u)\cup N(v)$. Let $R=T\backslash (N(u)\cup N(v))$ and $R_1, R_2, \cdots, R_r$ be the components of $G[R]$.

Note that $G[T]$ is connected. For each $i\in \{1, 2, \ldots, r\}$, $R_i$ has a neighbor, say $t_i$, in $N(u)\cup N(v)$. If $t_i$ is not complete to $R_i$, we may choose two adjacent vertices $x$ and $y$ in $R_i$ with $t_ix\in E(G)$ and $t_iy\notin E(G)$, then either $zut_ixy$ or $zvt_ixy$ is a $P_5$ of $G$. Therefore, $t_i$ must be complete to $R_i$, and so $G[R]$ is $K_3$-free to avoid a $K_1+(K_1\cup K_3)$.

Let $T_v=N_T(v)\setminus N_T(u)$. If $R$ is not anticomplete to $T_v$, let $r\in R$ and $t_v\in T_v$ be a pair of adjacent vertices, then $rt_vvzu$ is a $P_5$ in $G$. Thus $R$ is anticomplete to $T_v$, and consequently, $(N_T(u), R\cup T_v)$ is a partition of $T$ as desired.  This proves Lemma~\ref{P5K1K3K1freeC5}. \qed

\renewcommand{\baselinestretch}{1}
\begin{lemma}\label{P5K1K3K1freeantihole}
Suppose that $G$ is $C_5$-free and contains an odd antihole $A$ with at least seven vertices. Let $S$ be the set of vertices which are complete to $A$, and let $T=N(A)\backslash S$. Then $G[S]$ is $K_1\cup K_3$-free, $T$ can be partition into at most $2k+1$ independent sets, and $N^2(A)=\emptyset$.
\end{lemma}\renewcommand{\baselinestretch}{1.2}
\pf Suppose that $V(A)=\{v_1, v_2, \ldots, v_{2k+1}\}$, where $k\ge 3$ and $v_iv_{i+1}\not\in E(G)$ for each $i\in \{1, 2\ldots, 2k+1\}$. Since $G$ is $K_1+(K_1\cup K_3)$-free, it is certain that $G[S]$ is $K_1\cup K_3$-free.

Note that the vertex of $T$ is neither complete nor anticomplete to $A$. For each vertex $u$  of $T$, there must be an $i_u\in \{1, 2, \ldots, 2k+1\}$ such that  $uv_{i_u}\not\in E(G)$ and $uv_{i_u+1}\in E(G)$, and so $uv_{i_u+3}\in E(G)$ to avoid either a $C_5$ or a $P_5$ depending on whether $uv_{i_u+2}\in E(G)$ or not. For each $i\in \{1, 2, \ldots, 2k+1\}$, let
$$\mbox{$T_i=\{v|v\in T, vv_i\notin E(G)$ but $vv_{i+1}\in E(G)$ and $vv_{i+3}\in E(G)$\} }.$$
Thus $T=\cup_{1\le i\le 2k+1} T_i$. Since $G$ is $K_1+(K_1\cup K_3)$-free, each $T_i$ is independent,  otherwise $G[\{v_i, v_{i+1}, v_{i+3}, x, x'\}]=K_1+(K_1\cup K_3)$ for any two adjacent vertices $x$ and $x'$ of $T_i$. Hence $T$ can be partition into at most $2k+1$ independent sets.

\medskip

Suppose that $N^2(A)\neq \emptyset$. Let $v$ be a vertex in $N(A)$ that has a neighbor, say $x$, in $N^2(A)$. It is obvious that $v\notin S$, otherwise a $K_1+(K_1\cup K_3)$ appears in $G$. Without loss of generality, we suppose that $v\in T_1$. Thus either a $K_1+(K_1\cup K_3)$ appears on $\{v, v_2, v_4, v_{2k}, x\}$ whenever $vv_{2k}\in E(G)$,  or a $P_5=xvv_2v_{2k}v_1$ appears whenever $vv_{2k}\not\in E(G)$. Therefore, $N^2(A)=\emptyset$. This proves Lemma~\ref{P5K1K3K1freeantihole}. \qed

\medskip


We are ready to prove Theorem~\ref{main3}.

\noindent\textbf{Proof of Theorem~\ref{main3}.}  Let $G$ be a $\{P_5, K_1+(K_1\cup K_3\}$-free graph. We may suppose that $G$ is connected, contains no clique cut set, and is not perfect. Thus $G$ contains a 5-hole  or an odd antihole with at least 7 vertices.

First suppose that $G$ contains a 5-hole $C=v_1v_2v_3v_4v_5v_1$. Since $G$ is $P_5$-free, we have that $V(G)=V(C)\cup N(C)\cup N^2(C)\cup N^3(C)$. By Lemma~\ref{P5free}$(a)$, we have that
\begin{equation*}
\begin{aligned}
N(C) & =N_{\{1,2,3,4,5\}}(C)\bigcup\limits_{1\le i\le 5} (N_{\{i, i+2\}}(C)\cup N_{\{i, i+1, i+2\}}(C)\\
         &\ \ \ \cup N_{\{i, i+1, i+3\}}(C)\cup N_{\{i, i+1, i+2, i+3\}}(C)).
\end{aligned}
\end{equation*}
By Lemmas~\ref{P5K3}, \ref{P5K1K3}, and \ref{P5K1K3K1free}$(a)$, we have that $\chi(G[N_{\{1,2,3,4,5\}}(C)])\le 3(\omega(G)-3)$, \\ $\chi(G[\bigcup\limits_{1\le i\le 5} N_{\{i,i+2\}}(C)])\le 15$, and $\chi(G[\bigcup\limits_{1\le i\le 5}(N_{\{i,i+1,i+2\}}(C)\cup N_{\{i,i+1,i+3\}}(C)\cup N_{\{i,i+1,i+2,i+3\}}(C))])\le 5$. Therefore, $\chi(G[N(C)])\le 3\omega(G)+11$.

By Lemmas~\ref{P5K3} and \ref{P5K1K3K1freeC5}, $\chi(G[N^2(C)])\le 6$ and  $\chi(G[N^3(C)])\le 3$. Since $\bigcup\limits_{1\le i\le 5}N_{\{i,i+2\}}(C)$ is anticomplete to $N^2(C)\cup N^3(C)$ by Lemma~\ref{P5K1K3K1free}$(a)$, we can color the vertices of $V(C)\cup N^2(C)\cup N^3(C)$ with the 15 colors used on $\bigcup\limits_{1\le i\le 5}N_{\{i,i+2\}}(C)$.  Thus $\chi(G)\le 3\omega(G)+11$.

\medskip

Now we suppose that $G$ is $\{P_5, C_5, K_1+(K_1\cup K_3\}$-free and contains an odd antihole $A$   with $|A|=2k+1\ge 7$. Let $S\subseteq N(A)$ be the set of all vertices that are complete to $A$, and let $T=N(A)\backslash S$. By Lemma~\ref{P5K1K3K1freeantihole}, we have that $V(G)=A\cup N(A)$, $G[S]$ is $K_1\cup K_3$-free, and $T$ is the union of $2k+1$ independent sets.  Hence $\chi(G[S])\le 3(\omega(S)-1)\le 3(\omega(G)-k-1)$ by Lemma~\ref{P5K1K3}, and so $\chi(G)\le \chi(A)+ \chi(G[S])+\chi(G[T])\le (k+1)+3(\omega(G)-k-1)+(2k+1)< 3\omega(G)$. This proves Theorem~\ref{main3}. \qed

\end{document}